\documentclass[12pt]{elsarticle}

\usepackage{amssymb}

\usepackage{amssymb}
\usepackage{enumerate}
\usepackage{latexsym}
\usepackage[only,ninrm,elvrm,twlrm,sixrm,egtrm,tenrm]{rawfonts}
\usepackage{indentfirst}
\usepackage{amsmath,amsfonts}
\usepackage{algorithm}
\usepackage{algorithmic}
\usepackage{multirow}
\usepackage{graphicx,psfrag}
\usepackage{graphics}
\usepackage{makeidx}
\usepackage{ifpdf}
\usepackage{amsmath}
\usepackage{xcolor}
\usepackage{array}
\usepackage[all,cmtip,color]{xypic}
\usepackage{lipsum}

\newtheorem{thm}{Theorem}[section]

\newtheorem{lem}{Lemma}[section]

\def\qed{\nopagebreak\hfill{\rule{4pt}{7pt}}
\renewcommand{\thefootnote}{\fnsymbol{footnote}}
\renewcommand{\baselinestretch}{1.3}
\renewcommand{\thefootnote}{}
\medbreak}
\renewcommand{\baselinestretch}{1.1}
\def\pf{\noindent {\it Proof.} }

\begin{document}

\begin{frontmatter}

\title{On the complexity of $k$-rainbow cycle colouring problems}

\author{Shasha Li$^a$, Yongtang Shi$^{b,}$\footnote{The corresponding author.}, Jianhua Tu$^c$, Yan Zhao$^d$}

\address{\small $^a$ Ningbo Institute of
Technology, Zhejiang University, Ningbo 315100, China}
\address{\small $^b$ Center for Combinatorics and LPMC,
 Nankai University, Tianjin 300071, China}
\address{\small $^c$ School of Science, Beijing University of Chemical Technology,
Beijing 100029, China}
\address{\small $^d$ Department of Mathematics,
 Taizhou University, Taizhou 225300, China\\[2mm]
\small Email: lss@nit.zju.edu.cn; shi@nankai.edu.cn; tujh81@163.com; zhaoyan81.2008@163.com}

\begin{abstract}
An edge-coloured cycle is $rainbow$ if all edges of the cycle have
distinct colours. For $k\geq 1$, let $\mathcal{F}_{k}$ denote the
family of all graphs with the property that any $k$ vertices lie on
a cycle. For $G\in \mathcal{F}_{k}$, a $k$-$rainbow$ $cycle$
$colouring$ of $G$ is an edge-colouring such that any $k$ vertices
of $G$ lie on a rainbow cycle in $G$. The $k$-$rainbow$ $cycle$
$index$ of $G$, denoted by $crx_{k}(G)$, is the minimum number of
colours needed in a $k$-rainbow cycle colouring of $G$. In this
paper, we restrict our attention to the computational aspects of
$k$-rainbow cycle colouring. First, we prove that the problem of
deciding whether $crx_1=3$ can be solved in polynomial time, but
that of deciding whether $crx_1 \leq \ell$ is NP-Complete, where $\ell\geq
4$. Then we show that the problem of deciding whether $crx_2=3$ can
be solved in polynomial time, but those of deciding whether $crx_2
\leq 4$ or $5$ are NP-Complete. Furthermore, we also consider the
cases of $crx_3=3$ and $crx_3 \leq 4$. Finally, we prove that the
problem of deciding whether a given edge-colouring (with an
unbounded number of colours) of a graph is a $k$-rainbow cycle
colouring, is NP-Complete for $k=1$, $2$ and $3$, respectively.
Some open problems for further study are mentioned.\\[2mm]
\noindent{\bf Keywords:} rainbow cycle; $k$-rainbow cycle colouring; $k$-rainbow cycle index; polynomial time;
NP-Complete\\[2mm]
{\bf AMS Subject Classification 2010:} 05C15, 05C38, 68Q25.
\end{abstract}

\end{frontmatter}
\section{Introduction}
We follow the terminology and notations of \cite{Bondy} and all
graphs considered here are finite and simple.


A well-known result of Dirac \cite{Dirac} states
that for $k\geq 2$, any $k$ specified vertices in a $k$-connected graph are contained in a cycle.
Bondy and Lov\'asz \cite{BL1981} proved that for $k\geq 2$, any $k-1$ specified vertices in a $k$-connected non-bipartite
graph are contained in an odd cycle; and for $k\geq 3$, any $k$ specified
vertices in a $k$-connected graph are contained in an even cycle. Bollob\'as and Brightwell \cite{BB1993} showed the following result: if $n\geq k\geq 3$ and $d \geq 1$ are such that $s=\lceil\frac k {\lceil n/d\rceil-1}\rceil\geq 3$, then for any $k$ vertices of degree at least $d$ in a graph of order $n$, there exists a cycle containing at least $s$ of the vertices. Very recently, Liu \cite{Liu} studied the following problem: For $k\geq 1$, let $\mathcal{F}_k$ denote the family
of all graphs $G$ with the property that any $k$ vertices of $G$ belong to a cycle.
An edge-coloured cycle is $rainbow$ if all edges of the cycle have
distinct colours. Consider an
edge-colouring of $G\in \mathcal{F}_k$ such that, any $k$ vertices are contained in a rainbow cycle. What
is the minimum number of colours in such an edge-colouring?

For $G\in
\mathcal{F}_{k}$, a $k$-$rainbow$ $cycle$ $colouring$ of $G$ is an
edge-colouring such that any $k$ vertices of $G$ lie on a rainbow
cycle in $G$. A 1-rainbow cycle colouring is simply called a {\it
rainbow cycle colouring}. An edge-coloured graph $G$ is
$k$-$rainbow$ $cycle$ $connected$ if its colouring is a $k$-rainbow
cycle colouring. The $k$-$rainbow$ $cycle$ $index$ of $G$, denoted
by $crx_{k}(G)$, is the minimum number of colours needed in a
$k$-rainbow cycle colouring of $G$. Thus, $crx_{k}(G)$ is
well-defined if and only if $G\in \mathcal{F}_{k}$. In \cite{Liu}, Liu studied the $k$-rainbow
cycle index for some special classes of graphs.

This concept is related to the concept of rainbow connection number
of graphs, which was introduced by Chartrand et al. \cite{CJMZ}. An
edge-coloured graph is $rainbow$ if the colours of its edges are
distinct. An edge-coloured graph $G$ is $rainbow$ $connected$ if any
two vertices are connected by a rainbow path. In this case, the
colouring is called a $rainbow$ $colouring$ of $G$. The $rainbow$
$connection$ $number$ of a connected graph $G$, denoted by $rc(G)$,
is the minimum number of colours that are needed in order to make
$G$ rainbow connected. Later, Krivelevich and Yuster \cite{KY}
extend the concept of rainbow connection to the vertex version. For
more results on rainbow connection and rainbow vertex connection, we
refer to the survey \cite{LSS} and some recent papers for digraphs \cite{LLLS,LLMS}.

The computational complexity of the rainbow (vertex-) connection
number has been studied extensively. In \cite{Caro}, Caro et al.
conjectured that computing the rainbow connection number is an
NP-Hard problem, as well as that even deciding whether a graph has
rainbow connection number 2 is NP-Complete, which was confirmed by
Chakraborty et al. \cite{CFMY}. For the rainbow vertex-connection
number, Chen et al. \cite{LXS} showed that for a graph G, deciding
whether the rainbow vertex connection number equals to 2 is
NP-Complete. Actually, there are many other results on this topic,
we refer to the recent papers \cite{HLS,LLS,Lauri1} and the PhD thesis
of Juho Lauri \cite{Lauri}.

In this paper, we restrict our attention to the computational
aspects of $k$-rainbow cycle colouring of graphs. In Section $2$, we
prove that the problem of deciding whether $crx_1=3$ can be solved
in polynomial time, but that of deciding whether $crx_1 \leq \ell$ is
NP-Complete, where $\ell\geq 4$. In Section $3$, we show that the
problem of deciding whether $crx_2=3$ can be solved in polynomial
time, but those of deciding whether $crx_2 \leq 4$ and $5$ are
NP-Complete. In Section $4$, we show that it is easy to check
whether $crx_3=3$ and $crx_3 \leq 4$. In the last section, we turn
to the problem of deciding whether the given edge-colouring (with an
unbounded number of colours) of a graph is a $k$-rainbow cycle
colouring and we prove that the problem is NP-Complete for $k=1$,
$2$ and $3$.

\section{$1$-Rainbow cycle index}

In this section, we consider the problem of determining whether a
given graph $G$ has a $1$-rainbow cycle colouring with $\ell$ colours,
that is, determining whether $crx_{1}(G)\leq \ell$, where $\ell\geq 3$.

We first present a polynomial-time algorithm for the case $\ell=3$ of the above problem.

\noindent {\bf Algorithm:} Deciding Whether $crx_{1}(G)\leq 3$

\noindent {\bf INPUT:} a graph $G=(V,E)$

\noindent {\bf OUTPUT:} a $1$-rainbow cycle colouring function $c$ of $G$ with three colours $1,2$ and $3$
or the conclusion that $crx_{1}(G)> 4$

\noindent 1: $set\ \chi:=\emptyset,\ F:=\emptyset$

\noindent 2: ${\bf while}\ V\backslash F\ is\ nonempty\ {\bf do}$

\noindent 3: $\ \ pick\ a\ vertex\ v\ from\ V\backslash F$

\noindent 4: $\ \ {\bf if}\ v\ has\ two\ adjacent\ neighbours\ u\ and\  w\ {\bf then}$

\noindent 5: $\ \ \ \ {\bf if}\ uw\in \chi\ {\bf then}$

\noindent 6: $\ \ \ \ \ \ colour\ vu\ and\ vw\ such\ that\ \{c(vu),c(vw)\}=\{1,2,3\}\backslash \{c(uw)\}$

\noindent 7: $\ \ \ \ \ \ replace\ \chi\ by\ \chi\cup\{vu,vw\}$

\noindent 8: $\ \ \ \ {\bf else}$

\noindent 9: $\ \ \ \ \ \ set\ c(vu):=1,\ set\ c(vw):=2\ and\ set\ c(uw):=3$

\noindent 10: $\ \ \ \ \ \ replace\ \chi\ by\ \chi\cup\{vu,vw,uw\}$

\noindent 11: $\ \ \ \ {\bf end}\ {\bf if}$

\noindent 12: $\ \ \ \ replace\ F\ by\ F\cup\{v,u,w\}$

\noindent 13: $\ \ {\bf else}$

\noindent 14: $\ \ \ \ return\ (crx_1(G)>4)$

\noindent 15: $\ \ {\bf end}\ {\bf if}$

\noindent 16: ${\bf end}\ {\bf while}$

\noindent 17: ${\bf while}\ E\backslash \chi\ is\ nonempty\ {\bf do}$

\noindent 18: $\ \ pick\ an\ edge\ e\ from\ E\backslash \chi$

\noindent 19: $\ \ set\ c(e):=1$

\noindent 20: $\ \ replace\ \chi\ by\ \chi\cup\{e\}$

\noindent 21: ${\bf end}\ {\bf while}$

\noindent 22: $return\ (c)$

Note that, in the above algorithm, if a vertex $v\in F$, then there
is a rainbow triangle containing $v$. For step 5, if an edge $uw\in
\chi$, then $uw$ has been coloured and the vertices $u,w\in F$.
Conversely, if $v\notin F$, then the edges adjacent to $v$ do not
belong to $\chi$.

Moreover, the running time of the above algorithm is bounded by $\mathcal{O}(n^3m)$. Thus, we have the following theorem:

\begin{thm}\label{thm1}
Given a graph $G$, the problem of deciding whether $crx_1(G)\leq 3$ can be solved in polynomial time.
\end{thm}

For $\ell\geq 4$, the problem of deciding whether $crx_1(G)\leq \ell$ turns out to be NP-Complete.
We establish its NP-completeness by reducing the following problem to it.
{\bf The} ${\bf \ell}${\bf-vertex-colouring problem:} given a graph $G$ and an integer $\ell$, decide whether
there exists an assignment of at most $\ell$ colours to the vertices of $G$ such that no pair of
adjacent vertices are coloured the same, namely whether $\chi (G)\leq \ell$. It is known that this $\ell$-vertex-colouring problem
is NP-Complete for $\ell\geq 3$.

\begin{thm}\label{thm2}
Given a graph $G$ and an integer $\ell\geq 4$, the problem of deciding whether $crx_1(G)\leq \ell$ is NP-Complete.
\end{thm}
\pf Clearly the problem belongs to NP. Now given an instance $G=(V,E)$ of the $\ell$-vertex-colouring problem,
we construct a graph $G'=(V',E')$ such that $\chi (G)\leq \ell$ if and only if $crx_1(G')\leq \ell$.

We start by constructing a star, with one leaf vertex corresponding to every vertex $v\in V(G)$ and
an additional central vertex $a$. Then for every edge $v_iv_j\in E(G)$, add a $v_iv_j$-path of length $\ell-2$:
$v_iu_{i,j}^{1}u_{i,j}^{2}\ldots u_{i,j}^{\ell-3}v_j$. For every pair $v_i,v_j$ of vertices such that $v_iv_j\notin E(G)$,
add two $v_iv_j$-paths of length $2$: $v_iw_{i,j}^{1}v_j$ and $v_iw_{i,j}^{2}v_j$.

More formally, the vertex set $V'$ of $G'$ is defined as follows:
$$V'=V\cup \{a\}\cup U\cup W$$
$$U=\{u_{i,j}^{1},u_{i,j}^{2},\ldots, u_{i,j}^{\ell-3}:\ v_iv_j\in E(G)\ and\ i<j\}$$
$$W=\{w_{i,j}^{1},w_{i,j}^{2}:\ v_iv_j\notin E(G)\ and\ i<j\},$$
\noindent and the edge set $E'$ is defined as follows:
$$E'=E_{1}\cup E_{2}\cup E_{3}$$
$$E_{1}=\{av_{i}:\ v_{i}\in V\}$$
$$E_{2}=\{v_iu_{i,j}^{1},u_{i,j}^{t}u_{i,j}^{t+1},u_{i,j}^{\ell-3}v_j:\ 1\leq t\leq \ell-4,v_iv_j\in E(G)\ and\ i<j\}$$
$$E_{3}=\{v_iw_{i,j}^{1},w_{i,j}^{1}v_{j},v_iw_{i,j}^{2},w_{i,j}^{2}v_{j}:\ v_iv_j\notin E(G)\ and\ i<j\}.$$

Now, if $crx_1(G')\leq \ell$, let $c'$ be a $1$-rainbow cycle colouring of $G'$ using $\ell$ colours.
We define the vertex-colouring $c$ of $G$ by $c(v_i)=c'(av_i)$, for every $v_i\in V$.
Note that, for each edge $v_iv_j\in E(G)$ and the vertex $u_{i,j}^{1}\in U\subseteq V'$,
there is only one cycle of length at most $\ell$ containing the
vertex $u_{i,j}^{1}$, that is, $v_iu_{i,j}^{1}u_{i,j}^{2}\ldots u_{i,j}^{\ell-3}v_jav_{i}$.
So $c'(av_i)\neq c'(av_j)$, and thus $c(v_i)\neq c(v_j)$, that is, $c$ is a proper $\ell$-vertex-colouring.

In the other direction, assume that $\chi (G)\leq \ell$ and $c:V\rightarrow \{1,2,\ldots,\ell\}$
is a proper $\ell$-vertex-colouring of $G$. Define the edge-colouring $c'$ of $G'$ as follows:

\noindent $\bullet$ For each edge $av_{i}\in E_1$, we set $c'(av_i)=c(v_i)$.

\noindent $\bullet$ For any $v_iv_j\in E(G)$, $c(v_i)\neq c(v_j)$, and so $c'(av_i)\neq c'(av_j)$.
Except the colours $c(v_i)$ and $c(v_j)$, colour the $v_iv_j$-path
$v_iu_{i,j}^{1}\ldots u_{i,j}^{\ell-3}v_j$ with the remaining $(\ell-2)$ colours such that no two edges
on the path have the same colour, that is,
$\{c'(av_i),c'(v_iu_{i,j}^{1}),\ldots,c'(u_{i,j}^{\ell-3}v_j),c'(av_j)\}=\{1,2,\ldots,\ell\}$.

\noindent $\bullet$ Finally, for any $v_iv_j\notin E(G)$, set $c'(v_iw_{i,j}^{1})=1$,
$c'(v_iw_{i,j}^{2})=2$, $c'(w_{i,j}^{1}v_{j})=3$ and $c'(w_{i,j}^{2}v_{j})=4$.

It is easy to verify that this edge-colouring $c'$ is indeed a
$1$-rainbow cycle colouring of $G'$ using $\ell$ colours. This
completes the proof.\qed

\section{$2$-Rainbow cycle index}

Next, we consider the problem of determining whether $crx_{2}(G)\leq \ell$, where $\ell\geq 3$.

In \cite{Liu}, Liu proved the following result for complete
graphs, which will be used in the later proof. For completeness, we include the proof in \cite{Liu} here.

\begin{lem}[\cite{Liu}]\label{lem1}
For $n\geq 3$, $crx_1(K_{n})=crx_2(K_{n})=3$.
\end{lem}
\pf It suffices to show that $crx_2(K_n) \leq 3$ for each $n\geq 3$. We use induction on $n$ to show that there is a $2$-rainbow cycle colouring of $K_n$ using three colours. For $n=3$, we simply take the rainbow-coloured $K_3$. Now suppose $n\geq 4$ and let $u$ be a vertex of $G\cong K_n$. By induction, there exists a $2$-rainbow cycle colouring $c'$ for $G-u\cong K_{n-1}$, using three colours. Define a colouring $c$ of $G$ as follows. If $n-1$ is even, then take a perfect matching $M$ of $G-u$, and for $vw\in M$, let $c(uv)$ and $c(uw)$ be the two colours different from $c'(vw)$. If $n-1$ is odd, then take $x_1,x_2,x_3\in V(G-u)$ such that $c'(x_1x_2)$, $c'(x_2x_3)$ and $c'(x_3x_1)$ are distinct, and take the perfect matching $M'$ of $G-\{u,x_,x_2,x_3\}$. For the colouring $c$, colour the edges from $u$ to the vertices of $M'$ in the same way as before, and let $c(ux_i)=c'(x_{i+1}x_{i+2})$ for $i=1,2,3$ (indices are taken modulo 3). Clearly $c$ also uses three colours. Now by induction, any two vertices of $G-u$ lie in a
rainbow triangle. If $v$ is a vertex of $M$ ($n-1$ even) or $M'$ ($n- 1$ odd), then $u$ and $v$ lie in
the rainbow triangle formed by $u$ and the edge of $M$ or $M'$ incident with $v$. For $n- 1$ odd,
$ux_ix_{i+1}$ is a rainbow triangle for $i = 1, 2, 3$. Hence, $c$ is a $2$-rainbow cycle colouring of $K_n$, and we are done by induction.\qed

For $\ell=3$, we have the following result:

\begin{thm}\label{thm3}
For a given graph $G$, $crx_2(G)=3$ if and only if $G$ is a complete
graph of order at least 3.
\end{thm}
\pf If $G$ is a complete graph of order at least 3, by Lemma \ref{lem1} we know that $crx_2(G)=3$.

In the other direction, clearly $|V(G)|\geq 3$. Assume that $G$ is
not a complete graph, that is, there is a pair $(u,v)$ of vertices
such that $uv\notin E(G)$. Then the shortest cycle containing $u$
and $v$ has length at least $4$. Thus, we have $crx_2(G)\geq 4$, a
contradiction.\qed

Therefore, given a graph $G$, we can decide whether $crx_2(G)=3$ in
polynomial time. However, for $\ell=4$, we can show that the problem is
NP-Complete. We denote the set $\{1,2,\ldots,n\}$ by $[n]$.

\begin{thm}\label{thm4}
Given a graph $G$, the problem of deciding whether $crx_2(G)\leq 4$ is NP-Complete.
\end{thm}
\pf Clearly the problem is in NP. We prove the NP-Completeness by
reducing ``the $4$-vertex-colouring problem" to it. Let $G=(V,E)$ be
an instance of the $4$-vertex-colouring problem, where
$V=\{v_{1},v_{2},\ldots,v_{n}\}$. We construct a graph $G'=(V',E')$
such that $\chi (G)\leq 4$ if and only if $crx_2(G')\leq 4$.

The vertex set $V'$ of $G'$ is defined as follows:
$$V'=\overline{V}\cup U\cup W^1\cup W^2$$
$$\overline{V}=V\cup \{v_{n+1}\},\ U=\{u_{i,j}:\ v_iv_j\in E(G)\ and\ i<j\}\cup \{a\}$$
$$W^1=\{w_{i,j}^{1},w_{i,j}^{2}:\ v_iv_j\notin E(G)\ and\ i<j\},$$
$$W^2=\{w_{i,n+1}^{1},w_{i,n+1}^{2}:\ i\in[n+1]\}.$$
And the edge set $E'$ of $G'$ is defined as follows:
$$E'=E_{1}\cup E_{2}^1\cup E_{2}^2\cup E_{2}^3\cup E_{3}\cup E_{4}^1\cup E_{4}^2\cup E_{5}$$
$$E_{1}=\{xy:\ x,y\in U\}$$
$$E_{2}^1=\{xy:\ x,y\in W^1\},\ E_{2}^2=\{xy:\ x\in W^1,y\in W^2\}$$
$$E_{2}^3=\{w_{i,n+1}^{k}w_{j,n+1}^{l}:\ 1\leq i<j\leq n+1,1\leq k\leq 2,1\leq l\leq 2\}$$
$$E_{3}=\{v_iu_{i,j},u_{i,j}v_j:\ v_iv_j\in E(G)\ and\ i<j\}\cup \{av_{i}:\ i\in [n]\}$$
$$E_{4}^1=\{v_iw_{i,j}^{1},w_{i,j}^{1}v_{j},v_iw_{i,j}^{2},w_{i,j}^{2}v_{j}:\ v_iv_j\notin E(G)\ and\ i<j\}$$
$$E_{4}^2=\{v_iw_{i,n+1}^{1},v_{n+1}w_{i,n+1}^{1},v_iw_{i,n+1}^{2},v_{n+1}w_{i,n+1}^{2}:\ i\in [n]\}\cup$$
$$\{v_nw_{n+1,n+1}^{1},v_{n+1}w_{n+1,n+1}^{1},v_nw_{n+1,n+1}^{2},v_{n+1}w_{n+1,n+1}^{2}\}$$
$$E_{5}=\{xy:\ x\in U,y\in W^1\cup W^2\}.$$

Now, suppose $crx_2(G')\leq 4$ and $c'$ is a $2$-rainbow cycle colouring of $G'$ using $4$ colours.
We define the vertex-colouring $c$ of $G$ by $c(v_i)=c'(av_i)$, for every $v_i\in V$.
Note that, $G'[\overline{V}]$ is an empty graph and  moreover, for any $v_iv_j\in E(G)$, only
the vertices $u_{i,j}$ and $a$ are adjacent to both $v_i$ and $v_j$.
Therefore, there is only one cycle of length at most 4 containing both $v_i$ and $v_j$, namely $v_iu_{i,j}v_jav_i$.
So $c'(av_i)\neq c'(av_j)$, and thus $c(v_i)\neq c(v_j)$, that is, $c$ is a proper $4$-vertex-colouring.

In the other direction, assume that $\chi (G)\leq 4$ and $c:V\rightarrow \{c_1,c_2,c_3,c_4\}$
is a proper $4$-vertex-colouring of $G$. We define the edge-colouring $c'$ of $G'$ as follows:

\noindent $\bullet$ Since $G'[U]$ and $G'[W^1]$ are complete graphs, by Lemma \ref{lem1}, we can colour
$E_1$ and $E_{2}^1$ by $c_1,c_2$ and $c_3$ such that $c'(E_1)$ and $c'(E_{2}^1)$ are
$2$-rainbow cycle colourings of $G'[U]$ and $G'[W^1]$, respectively.

\noindent $\bullet$ For each edge $xy\in E_{2}^2$, where $x\in W^1$ and $y\in W^2$, if
$y=w_{i,n+1}^{1}$, set $c'(xy)=c_3$; otherwise, $y=w_{i,n+1}^{2}$ and set $c'(xy)=c_4$, for any $i\in \{1,\ldots,n,n+1\}$.

\noindent $\bullet$ For the edges in $E_{2}^3$, we set $c'(w_{i,n+1}^{1}w_{j,n+1}^{1})=c'(w_{i,n+1}^{2}w_{j,n+1}^{2})=c_3$,
and $c'(w_{i,n+1}^{1}w_{j,n+1}^{2})=c'(w_{i,n+1}^{2}w_{j,n+1}^{1})=c_4$, where $1\leq i<j\leq n+1$.

\noindent $\bullet$ For each edge $av_{i}\in E_3$, set $c'(av_i)=c(v_i)$.

\noindent $\ $ For any $v_iv_j\in E(G)$, since $c$ is a proper vertex-colouring, $c(v_i)\neq c(v_j)$,
and so $c'(av_i)\neq c'(av_j)$. Thus we can colour the edges $v_iu_{i,j}$ and $v_ju_{i,j}$ by
$\{c_1,c_2,c_3,c_4\}\backslash \{c'(av_i),c'(av_j)\}$ such that $v_iu_{i,j}v_jav_i$ is a rainbow cycle.

\noindent $\bullet$ For any $v_iv_j\notin E(G)$$(1\leq i<j\leq n)$, set $c'(v_iw_{i,j}^{1})=c_1$,
$c'(v_iw_{i,j}^{2})=c_2$, $c'(w_{i,j}^{1}v_{j})=c_3$ and $c'(w_{i,j}^{2}v_{j})=c_4$.

\noindent $\bullet$ For the edges in $E_{4}^{2}$, set $c'(v_iw_{i,n+1}^{1})=c_1$ and
$c'(v_iw_{i,n+1}^{2})=c_2$, for $1\leq i\leq n+1$; set $c'(v_{n+1}w_{i,n+1}^{1})=c_3$ and
$c'(v_{n+1}w_{i,n+1}^{2})=c_4$, for $1\leq i\leq n$; set $c'(v_{n}w_{n+1,n+1}^{1})=c_3$ and
$c'(v_{n}w_{n+1,n+1}^{2})=c_4$.

\noindent $\bullet$ Finally, for each edge $xy\in E_{5}$, where $x\in U$ and $y\in W^1\cup W^2$, if
$y=w_{i,j}^{1}$, set $c'(xy)=c_3$; otherwise, $y=w_{i,j}^{2}$ and set $c'(xy)=c_4$.

Next, we show that the $4$-edge-colouring $c'$ is a $2$-rainbow cycle colouring of $G'$.
Let $x,y$ be any two vertices of $G'$.

\textbf{Case 1:} $x,y\in \overline{V}$.

If $\{x,y\}=\{v_i,v_j\}$, where $v_iv_j\in E(G)$, then $v_iu_{i,j}v_jav_i$ is a rainbow cycle containing $x$ and $y$.

If $\{x,y\}=\{v_i,v_j\}$, where $v_iv_j\notin E(G)$, then $v_iw_{i,j}^1v_jw_{i,j}^2v_i$ is a rainbow cycle containing $x$ and $y$.

If $\{x,y\}=\{v_i,v_{n+1}\}$, where $1\leq i\leq n$, then $v_iw_{i,n+1}^1v_{n+1}w_{i,n+1}^2v_i$ is a rainbow cycle containing $x$ and $y$.

\textbf{Case 2:} $x,y\in U$ or $W^1$.

$c'(E_1)$ and $c'(E_{2}^1)$ are indeed $2$-rainbow cycle colourings of $G'[U]$ and $G'[W^1]$, respectively.

\textbf{Case 3:} $x,y\in W^2$ or $x\in W^1$ and $y\in W^2$.

If $\{x,y\}=\{w_{i,n+1}^1,w_{i,n+1}^2\}$ or $\{x,y\}=\{w_{i,n+1}^k,w_{j,n+1}^l\}$, where $1\leq i\neq j\leq n+1,1\leq k\leq 2$
and $1\leq l\leq 2$, then $w_{i,n+1}^1v_iw_{i,n+1}^2w_{j,n+1}^lw_{i,n+1}^1$ is a rainbow cycle containing $x$ and $y$.

If $x\in W^1$ and $y=w_{i,n+1}^{k}\in W^2$, where $1\leq i\leq n+1$ and $1\leq k\leq 2$,
then $xw_{i,n+1}^1v_iw_{i,n+1}^2x$ is a rainbow cycle containing $x$ and $y$.

\textbf{Case 4:} $x\in U$ and $y\in \overline{V}\cup W^1\cup W^2$.

If $y=v_i\in \overline{V}$ or $y=w_{i,n+1}^{k}\in W^2$, where $1\leq i\leq n+1$ and $1\leq k\leq 2$,
then $v_iw_{i,n+1}^1xw_{i,n+1}^2v_i$ is a rainbow cycle containing $x$ and $y$.

If $y=w_{i,j}^{k}\in W^1$, where $v_iv_j\notin E(G)$, $i<j$ and $1\leq k\leq 2$,
then $v_iw_{i,j}^1xw_{i,j}^2v_i$ is a rainbow cycle containing $x$ and $y$.

\textbf{Case 5:} $x\in \overline{V}$ and $y\in W^1\cup W^2$.

Let $x=v_i$$ (1\leq i\leq n+1)$.

If $y\in W^1$, then $v_iw_{i,n+1}^1yw_{i,n+1}^2v_i$ is a rainbow cycle containing $x$ and $y$.

If $y=w_{j,n+1}^{k}\in W^2$, where $1\leq j\leq n+1$, $i\neq j$ and $1\leq k\leq 2$,
then $v_iw_{i,n+1}^1w_{j,n+1}^{k}w_{i,n+1}^2v_i$ is a rainbow cycle containing $x$ and $y$.

If $y=w_{i,n+1}^{1}$ or $w_{i,n+1}^{2}\in W^2$, choose any one vertex $z$ from $W^1$ and
then $v_iw_{i,n+1}^1zw_{i,n+1}^2v_i$ is a rainbow cycle containing $x$ and $y$.

We have considered all the cases and so $c'$ is indeed a $2$-rainbow cycle colouring of $G'$.
The proof is complete.\qed

By means of a similar construction as that in the proof of Theorem \ref{thm4}, we
can obtain a polynomial reduction of ``the $5$-vertex-colouring
problem" to the following problem.

\begin{thm}\label{thm5}
Given a graph $G$, the problem of deciding whether $crx_2(G)\leq 5$ is NP-Complete.
\end{thm}
\pf Let $G=(V,E)$ be an instance of the $5$-vertex-colouring problem,
where $V=\{v_{1},v_{2},\ldots,v_{n}\}$. We can construct a graph $G'$ similar to the graph
in Theorem \ref{thm4}, except that
$U=\{u_{i,j}^1,u_{i,j}^2:\ v_iv_j\in E(G)\ and\ i<j\}\cup \{a\}$
and $E_{3}=\{v_iu_{i,j}^1,u_{i,j}^2v_j:\ v_iv_j\in E(G)\ and\ i<j\}\cup \{av_{i}:\ i\in [n]\}$.
Moreover, the edges in $E_1$ still form a complete graph $G'[U]$ and the edges in $E_5$ still form
a complete bipartite graph between the vertices in $U$ and $W^1\cup W^2$. The other vertex sets and
edge sets remain unchanged.

Note that, now $G'[\overline{V}]$ is still an empty graph and for any $v_iv_j\in E(G)$, only
one vertex $a$ is adjacent to both $v_i$ and $v_j$. Therefore,
any cycle of length at most 5 containing both $v_i$ and $v_j$ must contain the edges $av_i$ and $av_j$.
Thus, for any $2$-rainbow cycle colouring $c'$ of $G'$ using $5$ colours, $c'(av_i)\neq c'(av_j)$.
Similarly, define the vertex-colouring $c$ of $G$ by $c(v_i)=c'(av_i)$, for every $v_i\in V(G)$,
and $c$ is indeed a proper $5$-vertex-colouring.

In the other direction, let $c:V\rightarrow \{c_1,c_2,c_3,c_4,c_5\}$
is a proper $5$-vertex-colouring of $G$. The edge-colouring $c'$ of
$G'$ we will define is similar to that in the proof of Theorem \ref{thm4}.

In Theorem \ref{thm4}, we colour $E_1$ by $c_1,c_2$ and $c_3$ such that $c'(E_1)$ is
$2$-rainbow cycle colouring of $G'[U]$.
Now, for any $u_{i,j}^{k}u_{i',j'}^{l}\in E_1$ and $u_{i,j}^{k}a\in E_1$, where $\{i,j\}\neq \{i',j'\}$,
$1\leq k\leq 2$ and $1\leq l\leq 2$, let $c'(u_{i,j}^{k}u_{i',j'}^{l})$ and $c'(u_{i,j}^{k}a)$
be the colours of $u_{i,j}u_{i',j'}$ and $u_{i,j}a$ in Theorem \ref{thm4}, respectively.

For each edge $av_{i}\in E_3$, we still set $c'(av_i)=c(v_i)$.
For any $v_iv_j\in E(G)$, since $c(v_i)\neq c(v_j)$, $c'(av_i)\neq c'(av_j)$.
Now colour the edges $v_iu_{i,j}^1$, $u_{i,j}^1u_{i,j}^2$ and $v_ju_{i,j}^2$ by
$\{c_1,c_2,c_3,c_4,c_5\}\backslash \{c'(av_i),c'(av_j)\}$ such that $v_iu_{i,j}^1u_{i,j}^2v_jav_i$ is a rainbow cycle.

For each edge $xy\in E_{5}$, where $x\in U$ and $y\in W^1\cup W^2$, if
$y=w_{i,j}^{1}$, we still set $c'(xy)=c_3$; otherwise, $y=w_{i,j}^{2}$ and we still set $c'(xy)=c_4$.

The colours of the other edges remain unchanged.

Similarly, it can be verified that this edge-colouring $c'$ is indeed a $2$-rainbow cycle colouring
of $G'$ using five colours. The proof is complete.\qed

\section{$3$-Rainbow cycle index}

For a given graph $G$, it is easy to check whether $crx_3(G)=3$.

\begin{thm}\label{thm6}
For a given graph $G$, $crx_3(G)=3$ if and only if $G\in \{K_3,
K_4\}$.
\end{thm}
\pf Obviously $crx_3(K_3)=3$. For $K_4$, let $V(K_4)=\{v_1,v_2,v_3,v_4\}$.
Define an edge-colouring $c$ of $K_4$ as follows: $c(v_1v_2)=c(v_3v_4)=1$,
$c(v_1v_3)=c(v_2v_4)=2$ and $c(v_1v_4)=c(v_2v_3)=3$. It is easy to check that
$c$ is a $3$-rainbow cycle colouring of $K_4$, and so $crx_3(K_4)=3$.

Conversely, it is clear that if $crx_3(G)=3$, then $G$ is a complete
graph and $|V(G)|\geq 3$. Now, suppose that $G=K_n$, where $n\geq
5$. For any edge-colouring of $K_n$ with $3$ colours, since
$\delta=n-1\geq 4$, there always exists a triangle $S\subseteq
V(K_n)$, two edges of which have the same colour. Thus there is no
rainbow cycle of length $3$ containing $S$, and so $crx_3(K_n)\geq
4$, for $n\geq 5$. The proof is complete.\qed

Next, we consider the problem of deciding whether $crx_{3}(G)=4$.
Recall that for positive integers $t_i$, $1\leq i\leq k$, the
$Ramsey$ $number$ $r(t_1,t_2,\ldots,t_k)$ is the smallest integer
$n$ such that every $k$-edge-colouring $(E_1,E_2,\ldots,E_k)$ of
$K_n$ contains a complete subgraph on $t_i$ vertices all of whose
edges belong to $E_i$, for some $i$, $1\leq i\leq k$. In particular,
the Ramsey number $r_k=r(t_1,t_2,\ldots,t_k)$, where $t_i=3$, $1\leq
i\leq k$, satisfies the following inequality (\cite{Bondy}).

\begin{lem}\label{lem2}
$r_k\leq \lfloor k!e\rfloor+1$.
\end{lem}

Now, we give the following theorem.

\begin{thm}\label{thm7}
Given a graph $G$, if $|V(G)|\geq 66$, then $crx_3(G)>4$.
\end{thm}
\pf Let $|V(G)|=n$, and then $G\subseteq K_n$. By Lemma \ref{lem2}, $r_4\leq \lfloor 4!e\rfloor+1\leq 66$.
If $n\geq 66\geq r_4$, then for any edge-colouring of $K_n$ with $4$ colours, there always exists a
monochromatic copy of $K_3$, say with vertex set $S$. Thus any cycle of length at most $4$ containing $S$
must contain at least two edges with end-vertices in $S$, and so cannot be rainbow.
Hence, $crx_3(G)\geq crx_3(K_n)>4$. The proof is complete.\qed

Note that for $|V(G)|\leq 65$, we can use enumeration method and
check all the possible colourings of $E(G)$. Though the time is
bounded by a fixed integer $t$, it is possible that $t$ is too big
and unbearable. Thus, it is necessary to find a more effective
algorithm to decide whether $crx_3(G)=4$, for $|V(G)|\leq 65$.

For $\ell\geq 5$, the complexity of the problem of deciding whether $crx_3(G)\leq \ell$ remains unknown.

\section{Rainbow cycle colouring}

Suppose we are given an edge-colouring of the graph. Is it then easier to verify whether the colouring
is a $1$-rainbow cycle colouring? Clearly, if the number of colours is constant, then this problem
becomes easy, simply by means of an exhaustive search. However, if the colouring is arbitrary, the
problem becomes NP-Complete.

\begin{thm}\label{thm8}
Given an edge-coloured graph $G$, the problem of checking whether the given colouring is a $1$-rainbow cycle colouring
is NP-Complete.
\end{thm}

We establish its NP-Completeness by reducing the following problem
from \cite{CFMY} to it.

\begin{lem}\label{lem3}
The following problem is NP-Complete: Given an edge-coloured graph $G$ and two vertices $s,t$ of $G$,
decide whether there is a rainbow path connecting $s$ and $t$.
\end{lem}

\noindent $Proof$ $of\ Theorem\ \ref{thm8}.$ The problem clearly
belongs to NP. Now given a graph $G=(V,E)$ with two special vertices
$s$ and $t$ and an edge-colouring $c$ of $G$, we construct a graph
$G'=(V',E')$ and define an edge-colouring $c'$ of $G'$ such that $s$
and $t$ are rainbow connected in $G$ under $c$ if and only if $c'$
is a $1$-rainbow cycle colouring of $G'$.

Let $V=\{v_{1},v_{2},\ldots,v_{n-2},v_{n-1}=s,v_{n}=t\}$ be the vertex set of the original graph $G$.
We set $V'=V\cup U\cup W$, where $U=\{u_{1},u_{2},\ldots,u_{n-2}\}$ and
$W=\{w_{1},w_{2},\ldots,w_{2n-4}\}$, and
$E'=E\cup \{v_iu_i,u_iv_{i+1}:i\in [n-3]\}\cup \{v_{n-2}u_{n-2},u_{n-2}v_1\}\cup \{sw_1,w_{2n-4}t\}
\cup \{w_iw_{i+1}:i\in [2n-5]\}$.

The colouring $c'$ is defined as follows:

\noindent $\bullet$ All edges $e\in E$ retain the original colours, namely $c'(e)=c(e)$.

\noindent $\bullet$ The edges $sw_1$, $\{w_iw_{i+1}:i\in [2n-5]\}$ and $w_{2n-4}t$ are coloured with $(2n-3)$ new colours
$c'_1,c'_2,\ldots,c'_{2n-3}$, respectively.

\noindent $\bullet$ Finally, set $c'(v_iu_i)=c'_{2i-1}$ and $c'(u_iv_{i+1})=c'_{2i}$, for $i\in [n-3]$.
Set $c'(v_{n-2}u_{n-2})=c'_{2n-5}$ and $c'(u_{n-2}v_1)=c'_{2n-4}$.

Obviously, if there is a rainbow path $P$ from $s$ to $t$ in $G$ under $c$, then
$P\cup sw_1w_2w_3\ldots w_{2n-4}t$ is a rainbow cycle containing $\{s,t\}\cup W$.
Moreover, $v_1u_1v_2u_2v_3\ldots v_{n-3}u_{n-3}v_{n-2}u_{n-2}v_1$ is a rainbow cycle containing
$V\backslash \{s,t\}\cup U$. Thus, $c'$ is indeed a $1$-rainbow cycle colouring of $G'$.

Now, in turn, assume that $c'$ is a $1$-rainbow cycle colouring of $G'$. Then the vertex $w_1$
must lie on a rainbow cycle $l$ in $G'$. Obviously, any cycle containing $w_1$, including the cycle $l$,
must contain the path $sw_1w_2w_3\ldots w_{2n-4}t$, which uses up all the new colours:
$c'_1,c'_2,\ldots,c'_{2n-3}$. Thus $l\backslash \{w_{1},w_{2},\ldots,w_{2n-4}\}$ is actually a
rainbow path from $s$ to $t$ in $G$ under $c$.
The proof is complete.\qed

Intuitively, given an edge-coloured graph, if the colouring is arbitrary, the problem of deciding
whether the colouring is a $k$-rainbow cycle colouring $(k\geq 2)$ is not easier than that of deciding
whether the colouring is a $1$-rainbow cycle colouring.
Actually, we can show that, for $k=2$ and $3$, the problem is indeed NP-Complete.

Firstly, we prove the following claim.

\begin{lem}\label{lem4}
The first problem defined below is polynomially reducible to the second one:

\noindent {\bf Problem 1:} Given an edge-coloured graph $G$, decide whether any $k$ vertices of $G$
are connected by a rainbow path, where $k\geq 2$.

\noindent {\bf Problem 2:} Given an edge-coloured graph $G$, decide whether the given colouring is a $k$-rainbow cycle colouring,
where $k\geq 2$.
\end{lem}
\pf Given a graph $G=(V,E)$ and an edge-colouring $c$ of $G$,
we construct a graph $G'=(V',E')$ and define an edge-colouring $c'$ of $G'$ such that for the resulting edge-coloured graph
the answer for Problem $2$ is ``yes" if and only if the answer for Problem $1$ for the original edge-coloured graph
is ``yes".

Let $V(G)=\{v_{1},v_{2},\ldots,v_{n}\}$. Now we add two adjacent vertices $a$ and $b$ to $G$, and join $a$ and $b$ to
each vertex of $V(G)$, that is, $V'(G')=V\cup \{a,b\}$ and $E'(G')=E\cup \{av_i,bv_i:i\in [n]\}\cup \{ab\}$.

The colouring $c'$ is defined as follows:

\noindent $\bullet$ All edges $e\in E$ retain the original colours, namely $c'(e)=c(e)$.

\noindent $\bullet$ The edges $\{av_i:i\in [n]\}$ are coloured with a new colour $c'_{1}$.

\noindent $\bullet$ The edges $\{bv_i:i\in [n]\}$ are coloured with a new colour $c'_{2}$.

\noindent $\bullet$ The edge $ab$ is coloured with a new colour $c'_{3}$.

Firstly, assume that for any $k$ vertices of $V(G)$, there is a rainbow path $P$ connecting the $k$ vertices
in $G$. Let $x$ and $y$ be the ends of $P$. Then, obviously $xPyabx$ is a rainbow cycle in $G'$ containing the
vertices $a,b$ and the $k$ vertices. Thus, it is easy to see that $c'$ is indeed a $k$-rainbow cycle colouring of $G'$.

Conversely, assume that $c'$ is a $k$-rainbow cycle colouring of $G'$. Thus, for any $k$ vertices of $V(G)\subseteq V'(G')$,
there is a rainbow cycle $C$ containing the $k$ vertices. If the vertex $a$ belongs to $C$, since the edges
$av_1,av_2,\ldots,av_n$ have the same colour, the vertex $b$ must belong to $C$.
Then $C\backslash \{a,b\}$ is a rainbow path connecting the $k$ vertices in $G$.
If the vertices $a,b\notin V(C)$, then the deletion of any one edge from $C$ can result in a
rainbow path connecting the $k$ vertices in $G$.
The proof is complete.\qed

The case $k=2$ of Problem $1$ is exactly the problem in the following lemma, and its NP-Completeness has
been confirmed in \cite{CFMY}.

\begin{lem}\label{lem5}
The following problem is NP-Complete: Given an edge-coloured graph $G$,
check whether the given colouring makes $G$ rainbow connected.
\end{lem}

Obviously, Problem $2$ is in NP. Then Lemma \ref{lem5}, combined
with Lemma \ref{lem4}, yields the following theorem.
\begin{thm}\label{thm9}
Given an edge-coloured graph $G$, the problem of checking whether
the given colouring is a $2$-rainbow cycle colouring is NP-Complete.
\end{thm}

Next, we establish the NP-Completeness of the case $k=3$ of Problem $1$ by reducing the problem in
Lemma \ref{lem3} to it.

\begin{lem}\label{lem6}
The following problem is NP-Complete: Given an edge-coloured graph $G$,
decide whether any three vertices of $G$ are connected by a rainbow path.
\end{lem}
\pf Given a graph $G=(V,E)$ with two special vertices $s$ and $t$ and an edge-colouring $c$ of $G$,
we construct a graph $G'=(V',E')$ and define an edge-colouring $c'$ of $G'$ such that
$s$ and $t$ are rainbow connected in $G$ under $c$ if and only if the colouring $c'$ makes
any three vertices of $G'$ connected by a rainbow path.

Let $V=\{v_{1},v_{2},\ldots,v_{n-2},v_{n-1}=s,v_{n}=t\}$ be the vertex set of the original graph $G$.
We set $V'=\overline{V}\cup U\cup W\cup \{\hat{s},\hat{t},s,t\}$, where
$\overline{V}=V\backslash \{s,t\}=\{v_{1},v_{2},\ldots,v_{n-2}\}$,
$U=\{u^{1}_i,u^{2}_i:i\in [n-2]\}$ and $W=\{w^{1}_i,w^{2}_i:i\in [n-2]\}$, and
$E'=E\cup \{v_iu^{1}_i,v_iu^{2}_i:i\in [n-2]\}\cup \{v_iw^{1}_i,v_iw^{2}_i:i\in [n-2]\}
\cup \{u^{a}_iu^{b}_j,w^{a}_iw^{b}_j:i,j\in [n-2],\ a,b\in \{1,2\}\}
\cup \{u^{1}_{1}\hat{s},u^{1}_{1}s,w^{1}_{1}\hat{s},w^{1}_{1}s,\hat{s}s,
u^{1}_{n-2}\hat{t},u^{1}_{n-2}t,w^{1}_{n-2}\hat{t},w^{1}_{n-2}t,\hat{t}t\}$.

The colouring $c'$ is defined as follows:

\noindent $\bullet$ All edges $e\in E$ retain the original colours, namely $c'(e)=c(e)$.

\noindent $\bullet$ The edges $\hat{s}s$ and $\{v_iw^{1}_{i}:i\in [n-2]\}$ are coloured with a new colour $c'_{1}$.

\noindent $\bullet$ The edges $\hat{t}t$ and $\{v_iu^{2}_{i}:i\in [n-2]\}$ are coloured with a new colour $c'_{2}$.

\noindent $\bullet$ The edges $\{u^{1}_{1}\hat{s},u^{1}_{1}s,w^{1}_{1}\hat{s},w^{1}_{1}s,
u^{1}_{n-2}\hat{t},u^{1}_{n-2}t,w^{1}_{n-2}\hat{t},w^{1}_{n-2}t\}$ and $\{v_iu^{1}_{i}:i\in [n-2]\}$
are coloured with a new colour $c'_{3}$.

\noindent $\bullet$ The edges $\{u^{1}_iu^{2}_{i}:i\in \{2,3,\ldots,n-2\}\}$ are coloured with a new colour $c'_{4}$.

\noindent $\bullet$ The edges $u^{1}_{1}u^{2}_{1}$ and $\{u^{a}_iu^{b}_{j}:1\leq i< j\leq n-2,\ a,b\in \{1,2\}\}$
are coloured with a new colour $c'_{5}$.

\noindent $\bullet$ The edges $\{v_iw^{2}_{i}:i\in [n-2]\}$ are coloured with a new colour $c'_{6}$.

\noindent $\bullet$ The edges $\{w^{1}_iw^{2}_{i}:i\in [n-2]\}$ are coloured with a new colour $c'_{7}$.

\noindent $\bullet$ The edges $\{w^{a}_iw^{b}_{j}:1\leq i< j\leq n-2,\ a,b\in \{1,2\}\}$
are coloured with a new colour $c'_{8}$.

Next, we always let $i,j,k\in [n-2]$ and $a,b\in \{1,2\}$.

Firstly, suppose that there is a rainbow path $P$ from $s$ to $t$ in $G$ under $c$.
Let $x,y$ and $z$ be any three vertices of $G'$ and $S=\{x,y,z\}$.
Now let us prove that $x,y$ and $z$ are connected by a rainbow path under $c'$.

\textbf{Case 1:} $S\subseteq U$ or $S\subseteq W$.

If $S=\{u^1_{i},u^1_{j},u^2_{k}\}\subseteq U$, where $i\neq j\neq k$, then
$u^1_{i}u^1_{j}v_jw^1_jw^2_kv_ku^2_{k}$ is a rainbow path connecting $S$.

If $S=\{u^1_{i},u^1_{j},u^1_{k}\}\subseteq U$, where $i\neq j\neq k$ and w.l.o.g., $k\neq 1$, then
$u^1_{i}u^1_{j}v_jw^1_jw^2_kv_ku^2_{k}u^1_{k}$ is a rainbow path connecting $S$.

If $S=\{u^1_{i},u^1_{j},u^2_{j}\}\subseteq U$, where $i\neq j$, then
$u^1_{i}u^1_{j}u^2_{j}$ is a rainbow path connecting $S$ (if $j=1$, let the path be $u^1_{i}u^1_{j}v_ju^2_{j}$).

The other subcases $S=\{u^1_{i},u^2_{j},u^2_{k}\}$, $S=\{u^2_{i},u^2_{j},u^2_{k}\}$ and
$S\subseteq W$ are similar.

\textbf{Case 2:} $S\subseteq \overline{V}$.

Let $S=\{v_{i},v_{j},v_{k}\}$, where $i\neq j\neq k$, and then
$v_iu^1_{i}u^2_{j}v_jw^1_jw^2_kv_k$ is a rainbow path connecting $S$.

\textbf{Case 3:} $S\subseteq \{\hat{s},\hat{t},s,t\}$.

Obviously, $\hat{s}sPt\hat{t}$ contains a rainbow path connecting $S$.

\textbf{Case 4:} $x,y\in U$ and $z\in \overline{V}$ or $W$ or $\{\hat{s},\hat{t},s,t\}$.

Let $\{x,y\}=\{u^a_{i},u^b_{j}\}$.

If $z=v_k\in \overline{V}$, then $u^a_{i}u^b_{j}v_jw^1_jw^2_kv_k$ is a rainbow path connecting $S$
(if $k=i$ or $j$, let the path be $u^b_{j}u^a_{i}v_i$ or $u^a_{i}u^b_{j}v_j$).

If $z\in W$, then $u^a_{i}u^b_{j}v_jw^1_jz$ is a rainbow path connecting $S$.

If $z\in \{\hat{s},s\}$ (or $\{\hat{t},t\}$) and $a=b=1$, then w.l.o.g., let $j\neq 1$ and
$u^1_{i}u^1_{j}u^2_{j}v_jw^1_jw^1_1z$ (or $u^1_{i}u^1_{j}u^2_{j}v_jw^1_jw^1_{n-2}z$) is a rainbow path connecting $S$.

If $z\in \{\hat{s},s\}$ (or $\{\hat{t},t\}$) and w.l.o.g., $b=2$,
then $u^a_{i}u^b_{j}v_jw^1_jw^1_1z$ (or $u^a_{i}u^b_{j}v_jw^1_jw^1_{n-2}z$) is a rainbow path connecting $S$.

The case $x,y\in W$ and $z\in V'\backslash W$ is similar.

\textbf{Case 5:} $x,y\in \overline{V}$ and $z\in U$ or $W$ or $\{\hat{s},\hat{t},s,t\}$.

Let $\{x,y\}=\{v_{i},v_{j}\}$, where $i\neq j$.

If $z\in U$ (or $W$), then $v_iw^1_{i}w^2_{j}v_ju^1_jz$ (or $v_iu^1_{i}u^2_{j}v_jw^1_jz$)
is a rainbow path connecting $S$.

If $z\in \{\hat{s},s\}$ (or $\{\hat{t},t\}$), then
$v_iw^1_{i}w^2_{j}v_ju^2_ju^1_1z$ (or $v_iw^1_{i}w^2_{j}v_ju^2_ju^1_{n-2}z$) is a rainbow path connecting $S$.

\textbf{Case 6:} $x,y\in \{\hat{s},\hat{t},s,t\}$ and $z\in U$ or $W$ or $\overline{V}$.

If $z\in U$ (or $W$), then $zu^1_{1}\hat{s}sPt\hat{t}$ (or $zw^1_{1}\hat{s}sPt\hat{t}$)
contains a rainbow path connecting $S$.

If $z=v_i\in \overline{V}$, then $v_iw^2_{i}w^1_{1}\hat{s}sPt\hat{t}$ contains a rainbow path connecting $S$.

\textbf{Case 7:} $|U\cap S|\leq 1$, $|W\cap S|\leq 1$, $|\overline{V}\cap S|\leq 1$ and
$|\{\hat{s},\hat{t},s,t\}\cap S|\leq 1$.

If $x=u^a_{i}\in U$, $y=v_j\in \overline{V}$ and $z\in W$, then $v_ju^b_{j}u^a_{i}v_iw^1_{i}z$
is a rainbow path connecting $S$, where $a\neq b$ (if $i=j$, let the path be $u^a_{i}v_iw^1_{i}z$).

If $x\in U$, $y\in \{\hat{s},s\}$ (or $\{\hat{t},t\}$) and $z=w^a_{k}\in W$, then
$yw^1_{1}w^a_{k}v_ku^2_{k}x$ (or $yw^1_{n-2}w^a_{k}v_ku^2_{k}x$) is a rainbow path connecting $S$.

If $x\in U$, $y\in \{\hat{s},s\}$ (or $\{\hat{t},t\}$) and $z=v_j\in \overline{V}$, then
$yw^1_{1}w^2_{j}v_ju^2_{j}x$ (or $yw^1_{n-2}w^2_{j}v_ju^2_{j}x$) is a rainbow path connecting $S$.

If $x\in \{\hat{s},s\}$ (or $\{\hat{t},t\}$), $y=v_j\in \overline{V}$ and $z=w^a_{k}\in W$, then
$xu^1_{1}u^2_{j}v_jw^1_{j}w^a_{k}$ (or $xu^1_{n-2}u^2_{j}v_jw^1_{j}w^a_{k}$) is a rainbow path connecting $S$.

In a word, $S$ is always connected by a rainbow path.

Conversely, assume that the colouring $c'$ makes any three vertices of $G'$ connected by a rainbow path.
Thus, for $S=\{u^1_1,\hat{s},\hat{t}\}$, there is a rainbow path $P$ connecting $S$.
Note that the edges adjacent to $u^1_{1}$ are coloured either $c'_3$ or $c'_5$,
the edges connecting $U$ and $V'\backslash U$ are coloured either $c'_2$ or $c'_3$,
and the edges adjacent to $\hat{t}$ are coloured either $c'_2$ or $c'_3$.
It follows that $u^1_{1}$ must be one end of the path $P$. If $\hat{s}$ is the other end and
$\hat{t}$ is an internal vertex of $P$, then it is easy to check that $P=\hat{s}sP't\hat{t}u^1_{n-2}u^1_{1}$,
where $P'$ is exactly a rainbow path from $s$ to $t$ in the
original graph $G$. Similarly, if $\hat{t}$ is the other end and
$\hat{s}$ is an internal vertex of $P$, then $P=u^1_{1}\hat{s}sP't\hat{t}$,
where $P'$ is a rainbow path from $s$ to $t$ in the
original graph $G$. The proof is complete.\qed

Now from Lemma \ref{lem6} and Lemma \ref{lem4}, we can get the
following theorem.
However, for $k\geq 4$, the complexity of Problem $1$ and $2$ remains unknown.

\begin{thm}\label{thm10}
Given an edge-coloured graph $G$, the problem of checking whether the given colouring is a $3$-rainbow cycle colouring
is NP-Complete.\qed
\end{thm}

\section{Future work}
The $k$-rainbow cycle
index, $crx_{k}(G)$ studied in the paper is a new topic and there are so many properties can be investigated.
Furthermore, it would be interesting to study the parameter $\mathcal{F}_k$. Some properties of $\mathcal{F}_k$ are characterized by Liu in \cite{Liu}.
One referee pointed out to study the complexity for checking if a graph is in $\mathcal{F}_k$, for various $k$. Hope to come back in the future.

{\bf Acknowledgments.}
The authors would like to thank two anonymous referees for many helpful comments and suggestions. The authors also would like to thank Henry Liu for introducing to them the problem of the determination of the $k$-rainbow cycle index, and for his insightful comments and suggestions for improving the presentation of the paper.

Shasha Li was partially
supported by National Natural Science Foundation of China (No. 11301480),
Zhejiang Provincial Natural Science Foundation of China (No. LY18A010002),
and the Natural Science Foundation of Ningbo, China (No.2017A610132). Yongtang Shi was partially
supported by China--Slovenia bilateral project ``Some topics in modern graph theory"  (No. 12-6), the Natural Science
Foundation of Tianjin (No. 17JCQNJC00300) and the National Natural Science Foundation of China. J. Tu was partially supported by the National Natural Science Foundation of
China (No. 11201021), BUCT Fund for Disciplines Construction and
Development (Project No. 1524).  Yan Zhao was partially supported
by the Natural Science Foundation of Jiangsu Province(No.BK20160573),
and the Natural Science Foundation of the Jiangsu Higher Education Institutions of China(No.16KJD110005).

\end{document}